\documentclass[12pt]{article}
\usepackage{amsfonts, amssymb}

\def\FF{{\bf F}}
\def\GG{{\bf G}}
\def\mm{{\mathfrak m}}
\def\O{{\cal O}}
\def\PP{{\bf P}}
\def\Z{\mathbb Z}
\def\dim{\mathop{\rm dim}\nolimits}

\def\spec{\mathop{\rm Spec}\nolimits}
\let\hra\hookrightarrow
\let\ov\overline

\newtheorem{theorem}{Theorem}

\newtheorem{example}[theorem]{Example}
\newtheorem{lemma}[theorem]{Lemma}

\def\example{\refstepcounter{theorem}\paragraph{Example \thetheorem}}
\def\rem{\refstepcounter{theorem}\paragraph{Remark \thetheorem}}

\def\proof{\paragraph{Proof}}

\makeatletter \def\l@section{\@dottedtocline{1}{0em}{1.2em}} \makeatother

\begin{document}

\centerline{\huge Representability of $GL_E$}

\medskip

\centerline{\sl Nitin Nitsure} 

\bigskip

The main result of this note is the following theorem. 

\begin{theorem} \label{GL}
{\bf (Representability of the functor $GL_E$) }  
Let $S$ be a noetherian scheme, and $E$ a coherent $\O_S$-module. 
Let $GL_E$ denote the contrafunctor on $S$-schemes
which associates to any $S$-scheme $f:T\to S$ 
the group of all $\O_T$-linear
automorphisms of the pullback $E_T= f^*E$ 
(this functor is a sheaf in the fpqc topology).
Then $GL_E$ is representable by a group scheme over $S$ 
if and only if $E$ is locally free.
\end{theorem}

The `if' part is obvious. The main work
is in proving the `only if' part, for which 
we need various preliminaries.
The following lemma is standard, and is the first step in
the construction of a flattening stratification of a noetherian 
scheme $S$ for a coherent sheaf on $\PP^n_S$. 

\begin{lemma}\label{flattening} 
If $R$ is a noetherian local ring and $E$ a finite $R$-module, 
there exists an ideal $I\subset \mm$
with the following property: the module $E/IE$ is free over $R/I$, and 
for any ideal $J\subset R$, the module $E/JE$ is free over $R/J$ 
if and only if $I\subset J$. By its property, $I$ is unique.
\end{lemma}

\proof Define $I$ to be the ideal generated by the matrix entries 
of the map $\varphi: R^q\to R^p$ where 
$R^q\stackrel{\varphi}{\to}R^p \to E\to 0$ 
is an exact sequence in which $p$ is minimal (equal to
the dimension of the vector space $E/\mm E$ over $R/\mm$, where
$\mm$ denotes the maximal ideal in $R$). 
It can be seen that this $I$ has the desired property.
\hfill$\square$

\rem \label{freemodule} As a consequence, if $R$ is a 
noetherian local ring and $E$ a finite $R$-module such that  
$E/\mm^nE$ is a free module over $R/\mm^n$
for each $n\ge 2$, then $E$ is free over $R$.
For by the above lemma, $I \subset \mm^n$
for each $n\ge 2$, hence $I=0$.

\begin{lemma} \label{principal} {\rm (Srinivas)}
Let $R$ be an artin local ring with maximal ideal $\mm$, 
and let $E$ be a finite $R$-module, with corresponding ideal 
$I$ as in Lemma \ref{flattening}.
Suppose that the ideal $I$ is a principal ideal and $\mm I =0$. 
Then $E$ is isomorphic to a direct sum of the form 
$R^m\oplus (R/I)^n$, where $m,\, n$ are non-negative integers.
\end{lemma}

\proof Let $R^q\stackrel{\varphi}{\to}R^p \to E\to 0$
be an exact sequence of $R$-modules such that 
$p = \dim_{R/\mm}(E/\mm E)$. The ideal 
$I$ is generated by the matrix entries 
of the map $\varphi: R^q\to R^p$. By assumption,
there exists some $a\in \mm$ with $I=(a)$ and $\mm a=0$.
If $a=0$ then $E$ is free, so now assume $a\ne 0$.
Hence every non-zero element of $I$ is of the form $ua$ where
$u\in R-\mm$ is some unit of $R$.  
Hence the non-zero matrix entries of $\varphi: R^q\to R^p$ (if any)
are of the form $ua$. Hence there is another matrix
$\psi$ whose non-zero entries are units of $R$, with 
$\phi = a \psi$. Changing the free basis of $R^q$ and $R^p$
gives row and column operations on $\psi$, which can be used to 
put it in a block form 
$$\left(\begin{array}{ll}
1_{m \times m}      & 0_{m\times (q-m)}      \\
0_{(p-m) \times  m} & 0_{(p-m)\times (q-m)} 
\end{array}
\right)$$
The lemma follows. \hfill$\square$

\begin{lemma}\label{parabolic}
Let $S$ be a noetherian scheme, and let $E$ be a coherent $\O_S$-module.
Let $E'$ be a coherent subsheaf of $E$, such that
the quotient $E/E'$ is locally free. If $GL_E$ is representable, then the 
subfunctor $P$ of $GL_E$ which consists of automorphisms of $E$ 
(over base changes) which preserve $E'$ is also representable,
and is represented by a closed subgroup scheme of $GL_E$ over $S$.
\end{lemma}

\proof If $f: F' \to F$ is a homomorphism of coherent sheaves on 
a scheme $T$ such that $F$ is locally free, then $T$ has
a closed subscheme $T_0 \hra T$ with the universal property that 
$f$ vanishes identically under a base-change $T'\to T$ if and only
if it factors via $T_0 \hra T$. Applying this 
with $T= GL_E$, $F'= E'_T$, $F = (E/E')_T$,
and with $f: E'_T \to (E/E')_T$ the composite
$E'_T \to E_T \stackrel{u}{\to} E_T \to  (E/E')_T$
where $u:E_T  \to E_T $ is the universal family of automorphisms 
over $T=GL_E$, we get a closed subscheme $P\subset GL_E$
which has the desired properties. 
\hfill$\square$

\begin{lemma}\label{hence affine}
Let $X$ be a scheme, and $I\subset \O_X$ 
a quasi-coherent ideal sheaf, with $I^n=0$ for some $n\ge 1$. 
Suppose that the closed subscheme $Y\subset X$
defined by $I$ is affine. Then $X$ is affine.
\end{lemma}

\proof By induction on $n$, we can reduce to the case where $I^2=0$.
As $I^2=0$, $I$ becomes an $\O_Y$-module. As $I$ is quasi-coherent
over $\O_X$, it is quasi-coherent over $\O_Y$. If $F$ is any quasi-coherent
sheaf on $X$, then we have a short exact sequence
$0\to IF \to F \to F/IF \to 0$. As $I^2=0$, both
$IF$ and $F/IF$ are $\O_Y$-modules, and these are quasi-coherent.
Hence as $Y$ is affine,
$H^1(Y, IF) = H^1(Y, F/IF) = 0$. But these are just cohomologies over
the space $X$, as topologically 
$Y$ is $X$. Hence by the long exact sequence of
$0\to IF \to F \to F/IF \to 0$, it follows that 
$H^1(X,F)=0$. As this holds for every quasi-coherent 
$\O_X$-module, $X$ is affine by Serre's theorem. \hfill$\square$

\begin{lemma}\label{finitetype}
Let $A$ be a ring and $I\subset A$ an ideal with 
$I^n=0$ for some $n\ge 1$. Let $B$ be an $A$-algebra, such that
$B/IB$ is finite-type over $A$ (equivalently, over $A/I$). Let 
$b_1,\ldots,b_m\in B$ such that
$B/I = A[\ov{b}_1,\ldots,\ov{b}_m]$, where $\ov{b}_i\in B/I$
is the residue of $b_i$. Then $B$ is generated as an $A$-algebra 
by $b_1,\ldots,b_m$.
\end{lemma}

\proof By induction on $n$, we are reduced to the case where
$I^2=0$. As $B/I = A[\ov{b}_1,\ldots,\ov{b}_m]$, 
any $x\in B$ can be written as
$x = f(b_1,\ldots,b_m) + uy$ where $f$ is a polynomial
in $m$ variables over $A$, $u\in I$, and $y\in B$. Similarly,
$y = g(b_1,\ldots,b_m) + vz$ where $g$ is a polynomial
in $m$ variables over $A$, $v\in I$, and $z\in B$. As $I^2=0$,
we get $x = f(b_1,\ldots,b_m) +ug(b_1,\ldots,b_m)$.
Hence $B = A[b_1,\ldots,b_m]$. \hfill$\square$

\begin{lemma}\label{R/I is bad}
Let $R$ be an artin local ring with maximal ideal $\mm$, 
and let $0\ne I \subset \mm$ be a non-zero proper ideal. 
Let $E = (R/I)^n\oplus R^m$ where $n\ge 1$ and $m\ge 0$. 
Then the functor $GL_E$ is not representable.
\end{lemma}

\proof By Nakayama, $\mm I \ne I$, so
we can base-change to $R/\mm I$ and assume that
$\mm I=0$, in particular, $I^2 =0$. 
Suppose $GL_E$ is represented by 
a group-scheme $G$ over $R$. 
The restriction of $G$ to $R/I$ is the affine scheme
$GL_{n+m,\,R/I}$ over $R/I$, and $I$ is a nilpotent ideal, 
hence $G$ must be affine by Lemma \ref{hence affine}, and 
finite-type over $R$ by Lemma \ref{finitetype}. 
By Lemma \ref{parabolic}, the automorphisms which preserve $(R/I)^n\subset E$
are represented by a closed subgroup scheme $P\subset G$. 
Let $P =\spec(A)$ where $A$ is a finitely generated $R$-algebra.

The elements of the group $P(R)$ are matrices with the block form
$\left( \begin{array}{cc}
X & Y \\
0 & Z 
\end{array}
\right)$
where $X \in GL_n(R/I)$, $Y \in Hom(R^m, (R/I)^n) = (R/I)^{mn}$,
and $Z\in GL_m(R)$. 
Hence the elements $g\in P(R)$ which restrict to the identity 
in $P(R/I)$, that is, elements of the kernel of
$P(R)\to P(R/I)$, are exactly the elements of the form
$\left( \begin{array}{cc}
1 & 0 \\
0 & 1+ W 
\end{array}
\right)$
where $W \in M_m(I)$ is an arbitrary matrix with all entries in $I$.

The restriction of $P$ to
$R/I$ is the parabolic subgroup scheme $H \subset GL_{n+m, R/I}$
which preserves $(R/I)^n \subset (R/I)^{n+m}$, with
coordinate ring
$$B = R/I\,[x_{i,j},\, y_{i,\beta},\,  z_{\alpha,\beta}, \, 
\det(x_{i,j})^{-1},\,  \det(z_{\alpha,\beta})^{-1}]$$
where $1\le i,j \le n$, and $1\le \alpha, \beta \le m$.
As $B= A/IA$ where $I^2 =0$, by Lemma \ref{finitetype} we
get that 
$$A = R\,[x_{i,j},\, y_{i,\beta},\,  z_{\alpha,\beta}, \, 
\det(x_{i,j})^{-1},\,  \det(z_{\alpha,\beta})^{-1}]/J$$
for some ideal 
$J \subset I R\,[x_{i,j},\, y_{i,\beta},\,  z_{\alpha,\beta}, \, 
\det(x_{i,j})^{-1},\,  \det(z_{\alpha,\beta})^{-1}]$.
Let $V \in M_n(I)$ be any arbitrary $n\times n$-matrix over $I$.
We can define an $R$-algebra homomorphism
$A\to R$ by 
$$x_{i,j} \mapsto \delta_{i,j} + v_{i,j},~
y_{i,\beta} \mapsto 0, \mbox{ and }
z_{\alpha,\beta}\mapsto \delta_{\alpha,\beta}.$$
Modulo $I$, this specializes to identity, hence this contradicts
the above description of the kernel of $P(R) \to P(R/I)$.
This contradiction proves the lemma. \hfill$\square$

\bigskip

Now all the necessary preliminaries are in place for completing the proof
of the main result.

{\bf Proof of the Theorem \ref{GL}} 
Suppose that $E$ is not locally free. 
By first passing to the local ring of $S$ 
at some point where $E$ is not locally free and 
then going modulo a high power of the maximal ideal (see
Remark \ref{freemodule}), we can assume that
$S=\spec(R)$ where $R$ is an artin local ring, and $E$ is a finite 
$R$ module which is not free. 
Let $0\ne I \subset \mm$ be the ideal defined by $E$ as in Lemma
\ref{flattening}, where $\mm$ is the maximal ideal of $R$. 
Let $I = (a_1,\ldots, a_r)$ where $r$ is the smallest number of
generators needed to generate the ideal $I$. If $r\ge 2$,
let $J= (a_1,\ldots, a_{r-1}) \subset I$. Then going modulo
$J$ (that is, by base-changing to $R/J$), we are reduced to the 
case where $I$ is a principal ideal. By further going modulo
$\mm I$, we can assume $\mm I=0$. Hence by Lemma \ref{principal}, 
$E$ splits as a direct sum
$R^m\oplus (R/I)^n$, where $n \ge 1$ as $E$ is not free. 
Hence $GL_E$ is not representable by Lemma \ref{R/I is bad},
which completes the proof of the theorem. \hfill$\square$ 

\example The functor on commutative rings, defined by
$R \mapsto (R/2R)^{\times}$ (the multiplicative group of
units in the ring $R/2R$), is not representable by a scheme. 
This follows by taking $S= \spec(\Z)$ and $E = \Z/2\Z$ in
the Theorem \ref{GL}. A shorter direct proof is also possible 
in this example, by using discrete valuation rings instead 
of artin local rings.

{\bf Direct proof } 
If a group scheme
$G\to \spec(\Z)$ represents this functor, then the fiber
of $G$ over the closed point $(2)$ will be $\GG_{m,\FF_2}$,
while over the open complement $\spec(\Z) - (2)$, 
the restriction of $G$ will be 
trivial. Let $U$ be an affine open
neighbourhood in $G$ of the identity point 
$1 \in \GG_{m,\FF_2} \subset G$,
and let $x \in \GG_{m,\FF_2}$ be a closed point other than
$1$ which is in $U$ 
(the purpose of using an affine open $U$ is 
to avoid any assumption about separatedness of $G$).
The residue field $\kappa(x)$ at $x$ is a finite extension of 
$\FF_2$, hence separable over $\FF_2$. 
Let $A$ be the henselization of the local ring $\Z_{(2)}$ 
with respect to the residue field extension $\FF_2 \subset \kappa(x)$.
This is a discrete valuation ring of characteristic zero, 
with maximal ideal $2A$ as $A$ is \'etale over $\Z_{(2)}$, 
and residue field $\kappa(x)$.
Therefore, $G(\kappa(x)) = \kappa(x)^{\times} = 
(A/2A)^{\times} = G(A)$, and so $x$ uniquely prolongs to
an $A$-valued point of $G$, which we denote by $x'$. 
Note that $x': \spec A \to G$ factors through $U\subset G$. 
Therefore we have points $1$ and $x'$ of $U(A)$
which coincide over the generic point of $A$, but differ
over the special point. This contradicts the separatedness
of $U\to \spec(\Z)$. \hfill$\square$

\medskip

{\bf Acknowledgement } I had useful discussions on various aspects
with S.M. Bhatwadekar, H\'el\`ene Esnault, D.S. Nagaraj, 
Kapil Paranjape, and V. Srinivas, resulting in crucial 
inputs in the proof.

\bigskip

{\it Address : School of Mathematics, Tata Institute of Fundamental Research,
Homi Bhabha Road, Mumbai 400 005, India. e-mail:} 
{\tt nitsure@math.tifr.res.in}

\centerline{3 April 2002, 21 June 2002}

\end{document}